\documentclass{article}
\usepackage{amssymb,amsmath}
\usepackage{anysize}
\usepackage{psfrag}
\newtheorem{thm}{Theorem}[section]
\newtheorem{def.}{Definition}[section]

\newtheorem{cor}{Corollary}[section]

\numberwithin{table}{section}

\begin{document}
\title{On the Maximum Number of Colors for Links}
\author{Slavik Jablan\\
        The Mathematical Institute\\
        Knez Mihailova 36\\
        P.O. Box 367, 11001, Belgrade\\
        Serbia\\
        \texttt{sjablan@gmail.com}\\
        and\\
        Louis H. Kauffman\\
        Department of Mathematics, Statistics and Computer Science\\
        University of Illinois at Chicago\\
        851 S. Morgan St., Chicago IL 60607-7045\\
        USA\\
        \texttt{kauffman@uic.edu}\\
        and\\
        Pedro Lopes\\
        Center for Mathematical Analysis, Geometry and Dynamical Systems\\
        Departament of Mathematics\\
        Instituto Superior T\'ecnico\\
        1049-001 Lisboa\\
        Portugal\\
        \texttt{pelopes@math.ist.utl.pt}\\
}
\date{February 16, 2013}
\maketitle

\bigbreak

\begin{abstract}
For each odd prime $p$, and for each non-split link admitting non-trivial $p$-colorings, we prove that the maximum number of Fox colors is $p$. We also prove  that we can assemble a non-trivial $p$-coloring with any number of colors, from the minimum to the maximum number of colors. Furthermore, for any rational link, we prove that there exists a non-trivial coloring of any Schubert Normal Form of it, modulo its determinant, which uses all colors available. If this determinant is an odd prime, then any non-trivial coloring of this Schubert Normal Form, modulo the determinant, uses all available colors. We prove also that the number of crossings in the Schubert Normal Form equals twice the determinant of the link minus $2$. Facts about torus links and their coloring abilities are also proved.
\end{abstract}

\bigbreak


\bigbreak


\bigbreak

\section{Introduction} \label{sect:intro}

\begin{figure}[!ht]
	\psfrag{a}{\Huge$a$}
	\psfrag{b}{\Huge$b$}
	\psfrag{c}{\Huge$c=2b-a$}
	\centerline{\scalebox{.35}{\includegraphics{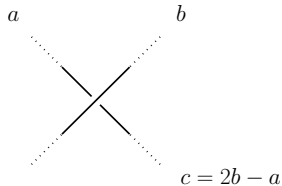}}}
	\caption{Detail of crossing of a diagram endowed with a non-trivial $p$-coloring, for odd prime $p$. The colors meeting at this crossing are all distinct.}\label{fig:x}
\end{figure}

Given a positive integer, say $p$, and a link $L$, a $p$-coloring of $L$ is an assignment of integers modulo $p$ to the arcs of a diagram of $L$ such that at each crossing, the integers assigned to the arcs meeting at this crossing satisfy the coloring condition. This condition states that, in the given modulus, twice the integer assigned to the over-arc equals the sum of the integers assigned to the under-arcs (see Figure \ref{fig:x}). In a coloring, the integers assigned to the arcs are called colors. Colorings were introduced by Fox in \cite{CFox} (see also \cite{SJablan} and \cite{lhKauffman}). The number of colorings is invariant under the Reidemeister moves and for any positive integer $p$ and for any link $L$ there are always the so-called trivial colorings. In a trivial coloring  each of the arcs of the diagram bears the same color. Furthermore, the fact that a link admits or not non-trivial colorings in a given modulus is also an invariant.
\noindent
In this sense, the interesting colorings are the non-trivial colorings i.e., those that use more than one color. How many distinct colors are needed to obtain a non-trivial coloring in a given modulus, for a given link? This question was first raised in \cite{Frank}. This number, called the minimum number of colors, is clearly a knot invariant. It is also a fascinating invariant from the point of view of its computation. In fact, by definition,  one has to find the minimum number of colors it takes to obtain a non-trivial coloring, for any given diagram of the link under study. This operation is to be repeated for each diagram of the link. Finally the answer is obtained by taking the minimum of all these minima. But any link has infinitely many diagrams. Therefore it is unreasonable to directly apply the definition to calculate this invariant, making this a challenging task.

\bigbreak

In \cite{kl} the last two authors obtained upper bounds for the minimum number of colors it takes to assemble a non-trivial coloring for the torus knots of type $(2, n)$. These upper bounds were improved in  \cite{klgame}. There is now a considerable number of references in the literature devoted to the minimum number of colors. Some are concerned with the calculation of (upper bounds of) minima for specific families of knots (\cite{kl}, \cite{lm2}). Others are concerned with establishing the exact value of the minima for specific modulus but for any link (\cite{kl}, \cite{satoh}, \cite{Oshiro}, \cite{Saito}, \cite{lm}). Others still are concerned with establishing transformations which allow the systematic reduction of the number of colors (\cite{kl}, \cite{klgame}).

\bigbreak

In this article we look into the dual problem. Given a positive modulus and a link, what is the maximum number of colors one can assign to the arcs of a diagram in order to obtain  a non-trivial coloring in the given modulus?

\begin{def.}[Maximum number of colors]
Given an integer $p>1$, assume there are non-trivial $p$-colorings
on a given link $L$. Assume further that $D$ is a diagram of $L$.
We let $n\sb{p, L}(D)$ denote the maximum number of distinct
colors assigned to the arcs of $D$ it takes to produce a
non-trivial $p$-coloring on $D$. We denote
$\mathrm{maxcol}\sb{p}L$ the maximum of these maxima over all
diagrams $D$ of $L$:
\[
\mathrm{maxcol}\sb{p}L := \max \{ n\sb{p, L}(D) \; | \; D
\textnormal{ is a diagram of } L \}
\]
For each $L$, we call $\mathrm{maxcol}\sb{p}L$ the {\bf maximum
number of colors of L, mod $p$}. In the sequel, we will drop the
``mod $p$'' whenever it is clear which $p$ we are referring to.
Note that $\mathrm{maxcol}\sb{p}L$ is tautologically a topological
invariant of $L$.
\end{def.}

On the other hand what are the possible numbers of colors a link can admit in order to put together a non-trivial coloring under a given modulus?

\begin{def.}[Spectrum of a Link for a given modulus]
Given an integer $p>1$, assume there are non-trivial $p$-colorings
on a given link $L$. We let
\[
\mathrm{Spec}\sb{p}L := \{ n\in \mathbf{Z}^+ \; | \;
\textnormal{there is a diagram, $D$, of $L$  with a non-trivial $p$-coloring which uses $n$ colors} \}
\]
For each $L$, we call $\mathrm{Spec}\sb{p}L$ the {\bf spectrum of L, mod $p$}. In the sequel, we will drop the
``mod $p$'' whenever it is clear which $p$ we are referring to.
Note that $\mathrm{Spec}\sb{p}L$ is a topological
invariant of $L$.
\end{def.}

\bigbreak

A {\bf split link} is a link which can be made to be contained in two disjoint balls in $3$-space. A {\bf non-split link} is a link which is not a split link.

For split links the  number of colorings in a given modulus equals the product of the number of colorings of each non-split component of the link. On the other hand,  for split links the question of minimum number of colors  is trivial. In fact, let us consider a split link which is then the  disjoint union of a certain number of non-split links. Coloring one of the split components with color $0$ and each of the others with color $1$ we obtain a non-trivial coloring of the link (in any modulus). We conclude that the minimum number of colors of any split link is $2$.

The question of the maximum number of colors for split links is also trivial thanks to a technique we develop below in the proof of Theorem \ref{thm:max}, see also Figures \ref{fig:x01}, \ref{fig:x01II}, and \ref{fig:x01IIn}. Again, on any given modulus and for any given split link, we color one of the split components  with color $0$ and all the others with color $1$. We pick two arcs bearing distinct colors within an appropriate neighborhood and perform a type II Reidemeister move by having the arc with color $1$ go over the arc with color $0$, thereby producing an arc with color $2$. Iterating this procedure we can give rise to new diagrams of this split link with any given number of distinct colors, up to the total number of colors available.

For these reasons,  in the present article we consider only non-split links.

\bigbreak
The questions that ultimately led to the results in this article were motivated by extensive calculations performed with the help of the programme \cite{LinKnot}. This programme includes functions that were written for the purpose of this article.

The main results of this article are the following.

\begin{thm}\label{thm:max}
Let $p$ be an odd prime. Let $L$ be a non-split link which admits non-trivial $p$-colorings.
Then
\[
\mathrm{maxcol}_p L = p
\]
\end{thm}

\bigbreak

\begin{cor}\label{cor:max}
We keep the hypothesis of Theorem \ref{thm:max} although allowing $p$ to be composite.
\begin{enumerate}
\item If there is a diagram of $L$ which admits a non-trivial $p$-coloring using the least number of colors and featuring  a crossing as in Figure \ref{fig:x01}, then the result follows.
\item If there is a diagram of $L$ which admits a non-trivial $p$-coloring using the least number of colors and featuring  a crossing as in Figure \ref{fig:x}, and $b-a$ is invertible mod $p$, then the result follows.
\end{enumerate}
\end{cor}

\bigbreak

Rational knots are also $2$-bridge links i.e., links admitting diagrams with exactly two over-arcs, see \cite{Schubert}. A Schubert Normal Form (SNF) is a specific sort of diagram of a $2$-bridge link with the following feature. When one travels along the diagram one goes under the two bridges, alternately. The SNF's of rational links  simplify the process of finding the diagram which maximizes the number of colors used. In the sequel we will be dealing exclusively with SNF whenever we are concerned with rational links. Also in the sequel ``rational link'' will stand for ``non-trivial rational link'' i.e., we will not consider the trivial knot nor the trivial $2$-component link when considering rational links.

\bigbreak

\begin{thm}\label{thm:maxrat}
Let $R$ be a rational link.

There is always a non-trivial coloring of an SNF of $R$, modulo its determinant, which uses all available colors.
\end{thm}

\bigbreak

\begin{cor}\label{cor:maxratprime}
Let $R$ be a rational knot with odd prime determinant $p$.

Any non-trivial $p$-coloring of an SNF of $R$ uses all available colors.
\end{cor}

\bigbreak

\begin{cor}\label{cor:N=2D-2}
Let $R$ be a rational link with determinant $D$. Let $N$ be the number of crossings in any of its SNF's. Then:
\[
N = 2D - 2
\]
\end{cor}

\bigbreak

Torus knots also provide interesting examples concerning the number of colors used in a coloring.  We remind the reader that
the classification of torus links according to their determinants has been done in \cite{Isidroetal} and that the braid closure of
\[
\bigl( \sigma_{r-1} \cdots \sigma_1\bigr)^{s}
\]
is our default diagram for the torus link $T(r, s)$.

\begin{thm}\label{thm:torus} Let $k$ and $l$ be positive integers.
\begin{enumerate}
\item Consider the diagram given by the braid closure of $\bigl( \sigma_{2l-1} \cdots \sigma_1\bigr)^{2k+1}$. It admits non-trivial $(2k+1)$-colorings using all colors available, with a uniform distribution of colors over the diagram.
\item Consider the diagram given by the braid closure of $\bigl( \sigma_{2k} \cdots \sigma_1\bigr)^{2l}$. It admits non-trivial $(2k+1)$-colorings using all colors available. For such a diagram, such a $(2k+1)$-coloring using all colors available does not exhibit a uniform distribution of the colors over this diagram. Furthermore, if additionally $\gcd(2l, 2k+1) =1$, a $(2k+1)$-coloring of such a diagram using all colors available cannot exhibit uniform distribution of the colors.
\item Consider the diagram given by the braid closure of $\bigl( \sigma_{2k-1} \cdots \sigma_1\bigr)^{2l}$. The determinant here is $2kl$. Modulo $2kl$ there is no uniform distribution of colors in a coloring using all colors available except when $k=1$.
\end{enumerate}
\end{thm}

\bigbreak

\begin{thm}\label{thm:spec}
Let $p$ be an odd prime. Let $L$ be a non-split link which admits non-trivial $p$-colorings.
Then
\[
\mathrm{Spec}_p L = \{ \mathrm{mincol}_p L , 1 + \mathrm{mincol}_p L , 2 + \mathrm{mincol}_p L , \dots , \mathrm{maxcol}_p L   \}
\]
\end{thm}

\bigbreak

\begin{cor}\label{cor:specbis}
We keep the hypothesis of Theorem \ref{thm:spec} although allowing $p$ to be composite.
\begin{enumerate}
\item If $L$ admits a minimal non-trivial $p$-coloring with a crossing as in Figure \ref{fig:x01}, then the result follows.
\item If $L$ admits a minimal non-trivial $p$-coloring with a crossing as in Figure \ref{fig:x}, and $b-a$ is invertible mod $p$, then the result follows.
\end{enumerate}
\end{cor}

\bigbreak

In Section \ref{sect:maxcolpL} we prove Theorem \ref{thm:max} and Corollary \ref{cor:max}; in Subsection \ref{subsect:funny} we look into what may happen if the hypothesis of Corollary \ref{cor:max} are not satisfied. In Section \ref{sect:maxrat} we prove Theorem \ref{thm:maxrat}, Corollary \ref{cor:maxratprime} and Corollary \ref{cor:N=2D-2}. In Section \ref{sect:torus} we prove Theorem \ref{thm:torus}. In Section \ref{sect:specpL} we prove Theorem \ref{thm:spec} and Corollary \ref{cor:specbis}. In Section \ref{sect:final} we collect a few questions for future work.

\bigbreak
\section{Proof of Theorem \ref{thm:max} and Corollary \ref{cor:max}.} \label{sect:maxcolpL}

\noindent

Proof (of Theorem \ref{thm:max}).

\bigbreak
We keep the notation and terminology of the statement of the Theorem.

If $L$ admits a non-trivial $p$-coloring, there is a diagram of $L$ along with assignments of integers mod $p$ to its arcs such that the coloring condition is satisfied at each crossing of the diagram and using more than one color in the process. In Figure \ref{fig:x} we depict a crossing of this diagram whose colors are all distinct.

\begin{figure}[!ht]
	\psfrag{a}{\Huge$0$}
	\psfrag{b}{\Huge$1$}
	\psfrag{c}{\Huge$2$}
	\psfrag{5}{\Huge$\mathbf{5}$}
	\centerline{\scalebox{.35}{\includegraphics{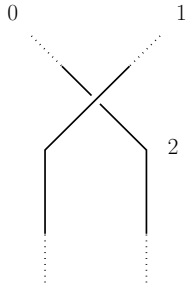}}}
	\caption{A new non-trivial $p$-coloring of the diagram. Now the crossing at issue bears the indicated colors, $0, 1, 2$.}\label{fig:x01}
\end{figure}

The $p$-colorings of a knot constitute a vector space since they are  solutions of a system of homogeneous linear equations over the integers mod $p$. Moreover, there are always the so-called trivial solutions, where each arc bears the same color. We now describe a method for passing from a non-trivial $p$-coloring with a crossing as depicted in Figure $\ref{fig:x}$ to another non-trivial $p$-coloring with the same crossing as depicted in Figure $\ref{fig:x01}$ i.e., $0\leftarrow a$, $1\leftarrow b$, and $2\leftarrow c$. This method is found in \cite{Saito} and \cite{lm}, we include it here for completeness.
\begin{enumerate}
\item If $a\neq 0$, subtract $a$ from each color in the diagram. If $a=0$ go to $2.$;
\item If $b\neq 1$, multiply each color in the coloring obtained in $1.$ by the inverse of $b$ mod $p$;
\end{enumerate}
After applying these two steps to the diagram under study the colors assigned to the crossing in Figure \ref{fig:x} become those depicted in   Figure \ref{fig:x01}.

\begin{figure}[!ht]
	\psfrag{a}{\Huge$0$}
	\psfrag{b}{\Huge$1$}
	\psfrag{c}{\Huge$2$}
	\psfrag{d}{\Huge$3$}
	\psfrag{5}{\Huge$\mathbf{5}$}
	\centerline{\scalebox{.35}{\includegraphics{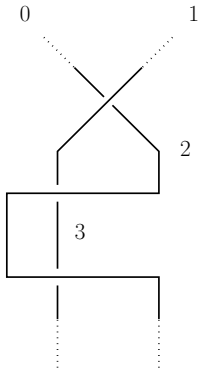}}}
	\caption{A new colored diagram only differing  from the one in Figure \ref{fig:x01} by a type II Reidemeister move and the insertion of color $3$. The present coloring includes definitely color $3$.}\label{fig:x01II}
\end{figure}

Suppose the diagram under study does not yet include color $3$, mod $p$; or it does include the color $3$ but locally does not look like Figure $3$. Then from this diagram we obtain a new diagram endowed with a non-trivial $p$-coloring and including color $3$ in such a way that around the crossing we have been addressing (see Figure \ref{fig:x01}), the new diagram now looks like Figure \ref{fig:x01II}.

\begin{figure}[!ht]
	\psfrag{b}{\Huge$k$}
	\psfrag{c}{\Huge$k+1$}
	\psfrag{d}{\Huge$k+2$}
	\psfrag{e}{\Huge$k+3$}
	\psfrag{5}{\Huge$\mathbf{5}$}
	\centerline{\scalebox{.35}{\includegraphics{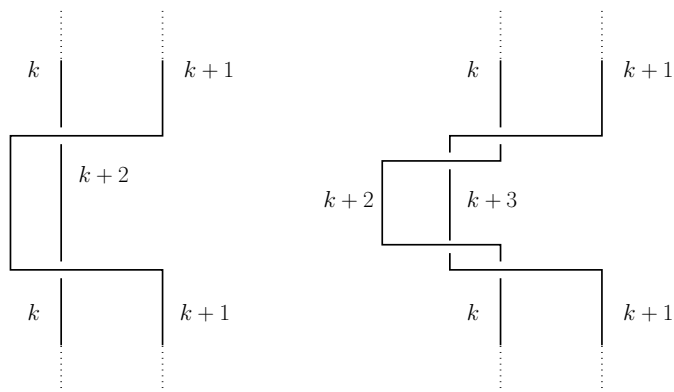}}}
	\caption{Obtaining a diagram involving colors $k, k+1, k+2, k+3$ (right-hand side) from a diagram involving colors $k, k+1, k+2$ (left-hand side).}\label{fig:x01IIn}
\end{figure}

We now prove that we can obtain inductively a diagram including all colors $0, 1, 2, \cdots , p-1$. In fact, if the current diagram involves a crossing including colors $k, k+1, k+2$ as depicted on the left-hand side of Figure \ref{fig:x01IIn}, then a type II Reidemeister move, along with the local change of colors so that the initial and final colorings are compliant, adds color $k+3$ to the existing colors, as the right-hand side of Figure \ref{fig:x01IIn} illustrates.
This completes the proof.
$\hfill \blacksquare$

\bigbreak

Proof (of Corollary \ref{cor:max}).

\bigbreak

Hypothesis $(1)$ in the statement of this Corollary ensures we have a crossing with colors as depicted in Figure \ref{fig:x01}.

Hypothesis $(2)$ implies that after subtracting each color of the non-trivial coloring at issue by $a$, and by multiplying each color of the resulting non-trivial coloring by the inverse of $b-a$ we obtain a diagram equipped with a crossing with colors as depicted in Figure \ref{fig:x01}.

From this point on we can apply Theorem \ref{thm:max} to prove this Corollary.
$\hfill \blacksquare$

\bigbreak
\subsection{What can go wrong when $b-a$ is not invertible mod $p$.} \label{subsect:funny}

\begin{figure}[!ht]
	\psfrag{0}{\Huge$0$}
	\psfrag{3}{\Huge$3$}
	\psfrag{6}{\Huge$6$}
	\psfrag{a}{\Huge$a$}
	\psfrag{b}{\Huge$b$}
	\psfrag{a+12(b-a)}{\Huge$a+12(b-a)$}
	\psfrag{b+12(b-a)}{\Huge$b+12(b-a)$}
	\psfrag{9}{\Huge$\mathbf{9}$}
	\centerline{\scalebox{.35}{\includegraphics{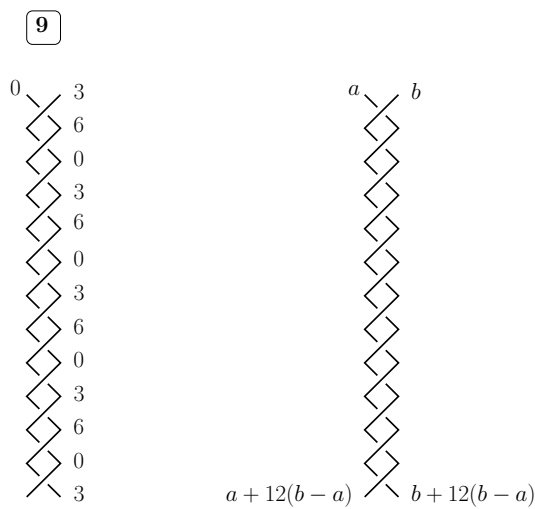}}}
	\caption{$T(2, 12)$ endowed with a non-trivial $9$-coloring involving $3$ colors (left-hand side). The propagation of colors down the $T(2, 12)$: $b-a$ has to be a multiple of $3$ for a non-trivial $9$-color to occur (right-hand side). Braid closure of each diagram is understood.}\label{fig:t12-9}
\end{figure}

\bigbreak

Inspection of Figure \ref{fig:t12-9} shows that, in the case of $T(2, 12)$ colored mod $9$, $b-a$ has to be a multiple of $3$, hence not invertible mod $9$. With such a choice for $b-a$ there are only $3$ colors in the coloring of $T(2, 12)$ ($0, 3, 6$, mod $9$, in the case at issue). The set $\{ 0, 3, 6\}$ is closed under the operation $a\ast b = 2b-a$ mod $9$. The same applies to the sets $\{ 1, 4, 7 \}$ or $\{ 2, 5, 8 \}$ mod $9$. Therefore, it is impossible to increase the number of colors starting from one of these colorings of this diagram.

\begin{figure}[!ht]
	\psfrag{0}{\Huge$0$}
	\psfrag{3}{\Huge$3$}
	\psfrag{6}{\Huge$6$}
	\psfrag{1}{\Huge$1$}
	\psfrag{2}{\Huge$2$}
	\psfrag{4}{\Huge$4$}
	\psfrag{5}{\Huge$5$}
	\psfrag{a}{\Huge$a$}
	\psfrag{b}{\Huge$b$}
	\psfrag{a+12(b-a)}{\Huge$a+12(b-a)$}
	\psfrag{b+12(b-a)}{\Huge$b+12(b-a)$}
	\psfrag{9}{\Huge$\mathbf{9}$}
	\centerline{\scalebox{.35}{\includegraphics{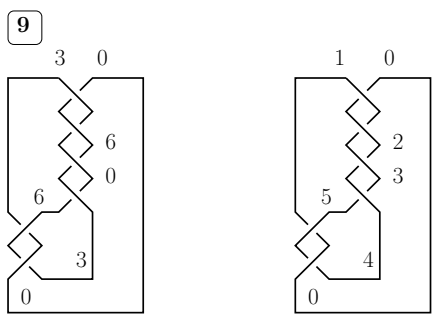}}}
	\caption{$R(4, 2)$ endowed with two essentially distinct non-trivial $9$-colorings.}\label{fig:r42-9}
\end{figure}

In Figure \ref{fig:r42-9} we see that, with $R(4, 2)$ we may increase the minimum number of colors, with a judicious choice of the coloring. Specifically, with the coloring on the left we are not able to increase the number of colors (arguing as above). For the coloring on the right it is possible to increase the number of colors using the technique set forth in the proof of Theorem \ref{thm:max}.

\bigbreak

\section{Proofs of Theorem \ref{thm:maxrat}, of Corollary \ref{cor:maxratprime}, and of Corollary \ref{cor:N=2D-2}.} \label{sect:maxrat}

\noindent

Proof (of Theorem \ref{thm:maxrat}).

\bigbreak
It is known that a rational link has either one or two components and that its coloring system of equations is equivalent to a single equation of the form
\[
D \times (b-a) = 0
\]
where $D$ is the determinant of the rational link and $a, b$ are variables corresponding to appropriate arcs of the diagram we are considering (see \cite{kldmtcs}, for example). We will now prove that we can choose these arcs to be the bridges when we represent a rational link by a Schubert Normal Form (SNF). These diagrams have the following feature. When traveling along the diagram, one goes under each of the bridges alternately (\cite{Schubert}).

\begin{figure}[!ht]
        \psfrag{br}{\Huge$b_r$}
        \psfrag{bl}{\Huge$b_l$}
        \psfrag{2br-bl}{\Huge$2b_r-b_l$}
        \psfrag{2bl-br}{\Huge$2b_l-b_r$}
        \psfrag{3br-2bl}{\Huge$3b_r-2b_l$}
        \psfrag{3bl-2br}{\Huge$3b_l-2b_r$}
        \psfrag{2kbr-(2k-1)bl}{\Huge$2kb_r-(2k-1)b_l$}
        \psfrag{2kbl-(2k-1)br}{\Huge$2kb_l-(2k-1)b_r$}
        \psfrag{(2k+1)br-2kbl}{\Huge$(2k+1)b_r-2kb_l$}
        \psfrag{(2k+1)bl-2kbr}{\Huge$(2k+1)b_l-2kb_r$}
        \centerline{\scalebox{.3}{\includegraphics{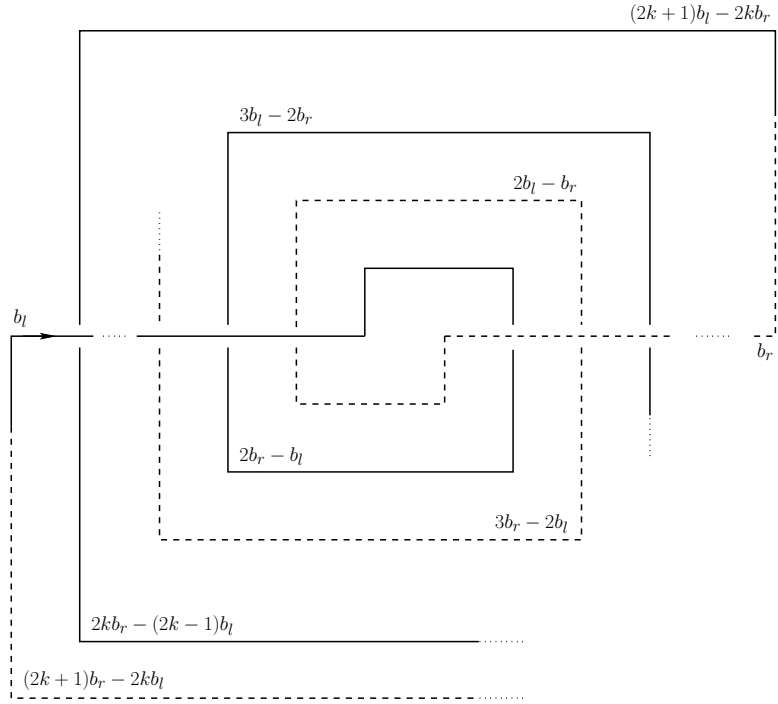}}}
        \caption{An SNF of a $1$-component rational knot. Its coloring system of equations is $(2k+1)(b_l-b_r) = 0$. The different line styles intend to partition the tour of the diagram into two halves. The one which begins at the bridge on the left- and the one which begins at the bridge on the right-hand side.}\label{fig:2k+1bridge}
\end{figure}

\begin{figure}[!ht]
        \psfrag{br}{\Huge$b_r$}
        \psfrag{bl}{\Huge$b_l$}
        \psfrag{2br-bl}{\Huge$2b_r-b_l$}
        \psfrag{2bl-br}{\Huge$2b_l-b_r$}
        \psfrag{3br-2bl}{\Huge$3b_r-2b_l$}
        \psfrag{3bl-2br}{\Huge$3b_l-2b_r$}
        \psfrag{2kbr-(2k-1)bl}{\Huge$2kb_r-(2k-1)b_l$}
        \psfrag{2kbl-(2k-1)br}{\Huge$2kb_l-(2k-1)b_r$}
        \centerline{\scalebox{.3}{\includegraphics{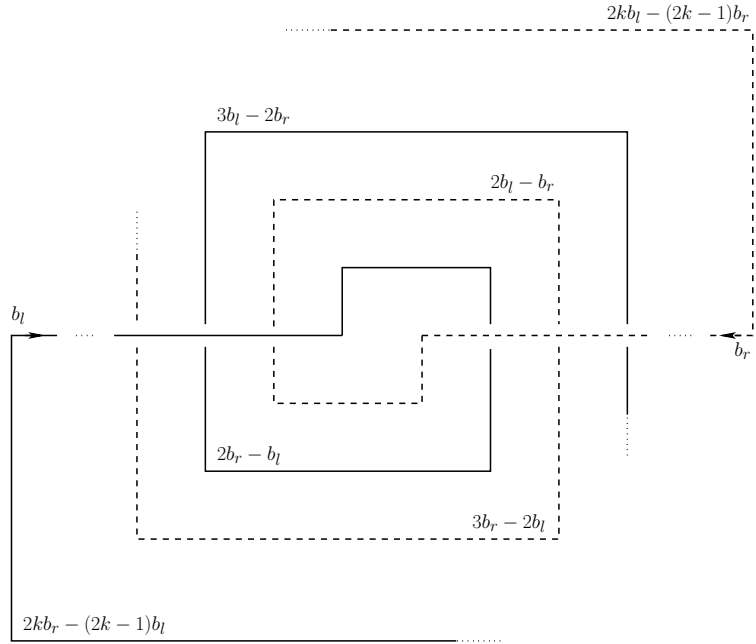}}}
        \caption{An SNF of a $2$-component rational link. Its coloring system of equations is $2k(b_l-b_r) = 0$. Different line styles denote different components of the link.}\label{fig:2kbridge}
\end{figure}

We will start by  describing how to retrieve the coloring system of equations from such a diagram. Let us assume we have a left bridge and a right bridge, depicted horizontally and that the variables assigned to them are $b_l$ and $b_r$, respectively. See Figures \ref{fig:2k+1bridge} and \ref{fig:2kbridge} for illustrative examples.

We start the tour of the diagram on the left bridge. Whenever we reach a crossing via an under-arc we write down the corresponding coloring equation, expressing the out-going under-arc in terms of the in-coming under-arc and the bridge. In particular, we do not mind crossings we arrive at via an over-arc. If the diagram at issue corresponds to a $2$-component link then, having started the tour on the left bridge we will end it on the left bridge; a complementary tour of the diagram (on the second component) will start and end on the right bridge (for instance). If the diagram corresponds to a $1$-component link, then having started on the left bridge the next over-arc we will reach is the right bridge. In either case we obtain two sets of equations as follows.

In the $1$-component case, the sequence of expressions associated to the under-arcs as we travel along the diagram is (see Figure \ref{fig:2k+1bridge} for an illustrative example):
\begin{align*}
&b_l\\
&2b_r-b_l \qquad \qquad \qquad \qquad (\text{$1$st time under right bridge})\\
&3b_l-2b_r \qquad \qquad \qquad \qquad (\text{$1$st time under left bridge})\\
&4b_r-3b_l \qquad \qquad \qquad \qquad (\text{$2$nd time under right bridge})\\
&5b_l-4b_r \qquad \qquad \qquad \qquad (\text{$2$nd time under left bridge})\\
& \qquad \qquad \qquad \qquad \dots \\
&2kb_r-(2k-1)b_l \qquad \qquad (\text{$k$-th time under right bridge})\\
&(2k+1)b_l-2kb_r \qquad \qquad (\text{$k$-th time under left bridge})
\end{align*}
so if the first arrival on the right bridge occurs after the $k$-th time under the left bridge we obtain the equation $(2k+1)b_l-2kb_r = b_r$ i.e.,
\begin{align*}
(2k+1)(b_l-b_r) = 0
\end{align*}

Now starting at the right bridge and doing the rest of the tour we obtain analogously
\begin{align*}
(2k'+1)(b_l-b_r) = 0
\end{align*}

Now in an SNF the under-segment which stems from the left bridge and the under-segment which stems from the right bridge are found alternately at crossings as one travels along the bridges of the diagram (\cite{Schubert}, page 141). Hence $k'=k$ and the system formed by the two equations above is equivalent to:

\[
(2k+1) (b_l-b_r) = 0
\]
We note that $2k$ is half the number of crossings in the SNF, in this instance.
\bigbreak

For the $2$-component case the analysis is similar and we obtain a coloring system of equations formed by
\begin{align*}
&2k(b_l-b_r) = 0\\
&2k'(b_l-b_r) = 0
\end{align*}
which is equivalent to
\begin{align*}
(2k)  (b_l-b_r) = 0
\end{align*}
since $k'=k$ arguing as before. We note that here $2k-1$ is half the number of crossings in the SNF.

\bigbreak

In either case we obtain a coloring system of equations made up of a single equation of the form $D (b_l - b_r) = 0$. Therefore $D$ is the determinant of the link and any $m$-coloring of the diagram is uniquely determined once we assign colors to $b_l$ and $b_r$ satisfying $D (b_l - b_r) = 0$ mod $m$. In passing, if $N$ is the number of crossings in the SNF then
\[
D = \frac{N}{2}+1
\]
in both instances, which concludes the proof of Corollary \ref{cor:N=2D-2}.
\bigbreak

Returning to the proof of Theorem \ref{thm:maxrat},
in the sequel we will then use $0$ and $1$ for the colors on the bridges of our $2$-bridge diagrams to obtain a non-trivial coloring, modulo the determinant, $D$, of the $2$-bridge link under study.

Consider a rational link with determinant $D$  represented  by a $2$-bridge diagram. As shown above, we obtain a non-trivial coloring modulo $D$ by  assigning $b_l =0, b_r=1$.

Let C be the set of colors (mod $D$) of the current non-trivial $D$-coloring of the $2$-bridge diagram of $R$. We remark that, with the exception of the bridges, each under-arc is an under-arc of both bridges.

Suppose color $i$ is missing from this diagram i.e., $i\notin C$ (note that $0\neq i \neq 1$).

\begin{figure}[!ht]
        \psfrag{0}{\Huge$0$}
        \psfrag{1}{\Huge$1$}
        \psfrag{i}{\Huge$i$}
        \psfrag{-i}{\Huge$-i$}
        \psfrag{-i+2}{\Huge$-i+2$}
        \centerline{\scalebox{.35}{\includegraphics{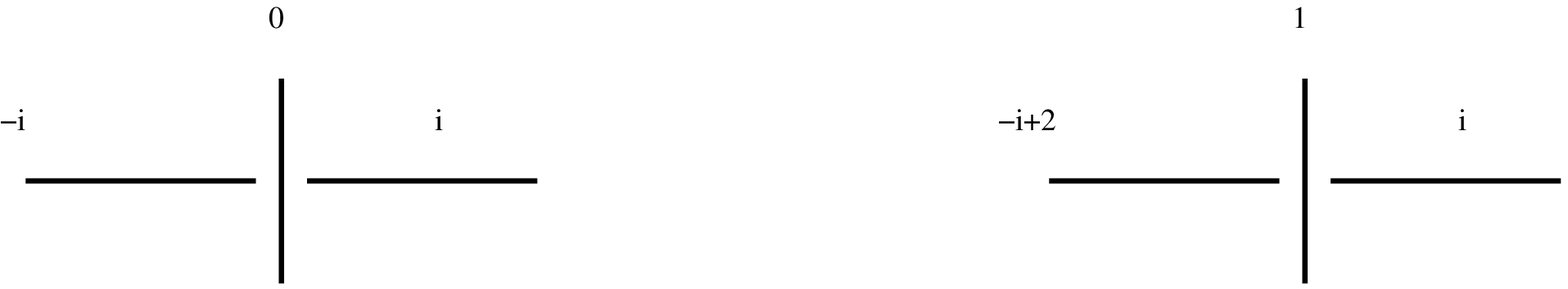}}}
        \caption{If $-i\in C$ then $i\in C$ (left); if $-i+2 \in C$ then $i \in C$ (right).}\label{fig:i1}
\end{figure}

\begin{figure}[!ht]
        \psfrag{0}{\Huge$0$}
        \psfrag{1}{\Huge$1$}
        \psfrag{i}{\Huge$i$}
        \psfrag{-i}{\Huge$i+2$}
        \psfrag{-i+2}{\Huge$i-2$}
        \psfrag{i+2}{\Huge$-i$}
        \psfrag{i-2}{\Huge$-i+2$}
        \centerline{\scalebox{.35}{\includegraphics{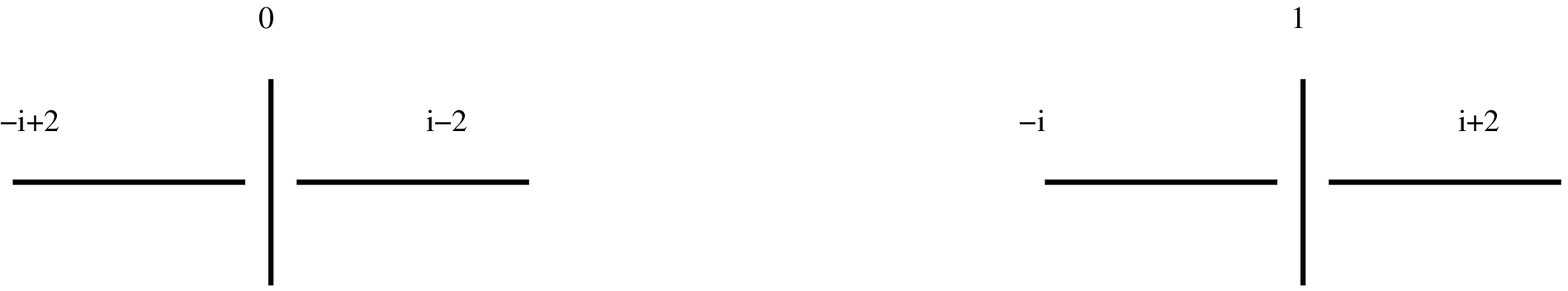}}}
        \caption{If $i - 2 \in C$ then $- i + 2 \in C$ (left); if $i+2 \in C$ then $-i \in C$ (right).}\label{fig:i2}
\end{figure}

Then (mod $D$) $-i$ and $2-i$ are not in $C$ for otherwise $i$ would also be in $C$ (see Figure \ref{fig:i1}).
But if $-i$ and $2-i$ are not in $C$, then $i+2$ and $i-2$ are also not in C (see Figure \ref{fig:i2}).
Therefore, by induction if $i \notin C$ then $i+2\notin C, i-2\notin C$; $i+4\notin C, i-4\notin C$, and so on and so forth. If $D$ is odd, then colors $0$ and $1$ will also not be present which is absurd. If $D$ is even at least one of them will not be present, which is also absurd. This completes the proof.

$\hfill \blacksquare$

Proof (of Corollary \ref{cor:maxratprime}).

\bigbreak
We keep the notation and terminology of the statement of the Corollary; in particular $p$ is an odd prime.

Consider a non-trivial $p$-coloring of an SNF of the rational knot $R$. We will first prove that the bridges cannot bear the same color. Assume to the contrary and without loss of generality, suppose this same color born by the bridges is $0$. Since the $p$-coloring is non-trivial then one of the under-arcs bears color $i\neq 0$. Starting from this under-arc and going under a bridge, the next under-arc is colored $-i\neq 0$. Hence, each under-arc is colored either $i\neq 0$ or $-i\neq 0$. Since the bridges are also under-arcs at some crossings then the bridges would also have to be colored $i\neq 0$ or $-i\neq 0$, which is absurd. Therefore, the bridges bear distinct colors when the diagram is endowed with a non-trivial $p$-coloring.

If the non-trivial $p$-coloring endows the bridges with colors $0$ and $1$, then we can use the proof of the Theorem to conclude that this coloring uses all colors available.

Consider now a non-trivial $p$-coloring where the bridges bear distinct colors $\{ b_l, b_r \} \neq \{ 0, 1 \}$ modulo $p$. We are going to map this $p$-coloring generated by the colors $b_l, b_r$ on the bridges, to a $p$-coloring generated by the colors $0, 1$ on the bridges, and prove that this mapping is bijective, thus proving that the original coloring generated by colors $b_l, b_r$ uses all colors available. Subtract each color of the coloring generated by bridges bearing $b_l, b_r$ by color $b_l$. Multiply each color of the coloring so obtained by $(b_r-b_l)^{-1}$ modulo $p$. The $p$-coloring now obtained is a non-trivial $p$-coloring with colors $0, 1$ on the bridges. Thus, it uses all colors available, as we saw before. To obtain the colors of the original coloring (the $c_i$'s), from the colors of this coloring (the $c'_i$'s), the expression is:
\[
c_i = (b_r-b_l)c'_i + b_l
\]
If $c_i = c_j$ then
\[
(b_r-b_l)c'_i + b_l = (b_r-b_l)c'_j + b_l \qquad \Leftrightarrow  \qquad (b_r-b_l)c'_i = (b_r-b_l)c'_j  \qquad \Leftrightarrow \qquad  c'_i = c'_j
\]
The proof is complete.

$\hfill \blacksquare$

\bigbreak

\section{Proof of Theorem \ref{thm:torus}.} \label{sect:torus}

\noindent

Proof (of Theorem \ref{thm:torus}).

\bigbreak
For this proof we will refer  to Figures \ref{fig:t2l2k+1} and \ref{fig:t2k+12l}. In these Figures the braid closure of the braids there depicted is understood. The boxed expressions on the top left of these Figures stand for the moduli with respect to which the colorings of the diagrams are being considered.

\begin{enumerate}
\item View Figure \ref{fig:t2l2k+1}. We believe this Figure clarifies the existence of a $(2k+1)$-coloring using all colors available, as stated. Now for the issue of color uniformity. Consider the bridge colored with color $s$, modulo $2k+1$. Color $s$ appears exactly $2l-1$ times on the diagram, $2(l-1)$ times right above and right below the indicated bridge and once on the bridge. Thus, each color appears exactly $2l-1$ times on the diagram.

    It has been brought to our attention that this non-trivial coloring is essentially the same as the one found on page 6 of the pre-print \cite{breil-oesp-taal}.
\begin{figure}[!ht]
	\psfrag{T2l2k+1}{\Huge$\mathbf{T(2l, 2k+1)}$}
	\psfrag{2k+1}{\Huge$\mathbf{2k+1}$}
	\psfrag{2k}{\Huge$2k$}
	\psfrag{2k-1}{\Huge$2k-1$}
	\psfrag{2k-2}{\Huge$2k-2$}
	\psfrag{0}{\Huge$0$}
	\psfrag{1}{\Huge$1$}
	\psfrag{2}{\Huge$2$}
	\psfrag{3}{\Huge$3$}
	\psfrag{4}{\Huge$4$}
	\psfrag{5}{\Huge$\mathbf{5}$}
	\psfrag{6}{\Huge$6$}
	\psfrag{7}{\Huge$\mathbf{7}$}
	\centerline{\scalebox{.25}{\includegraphics{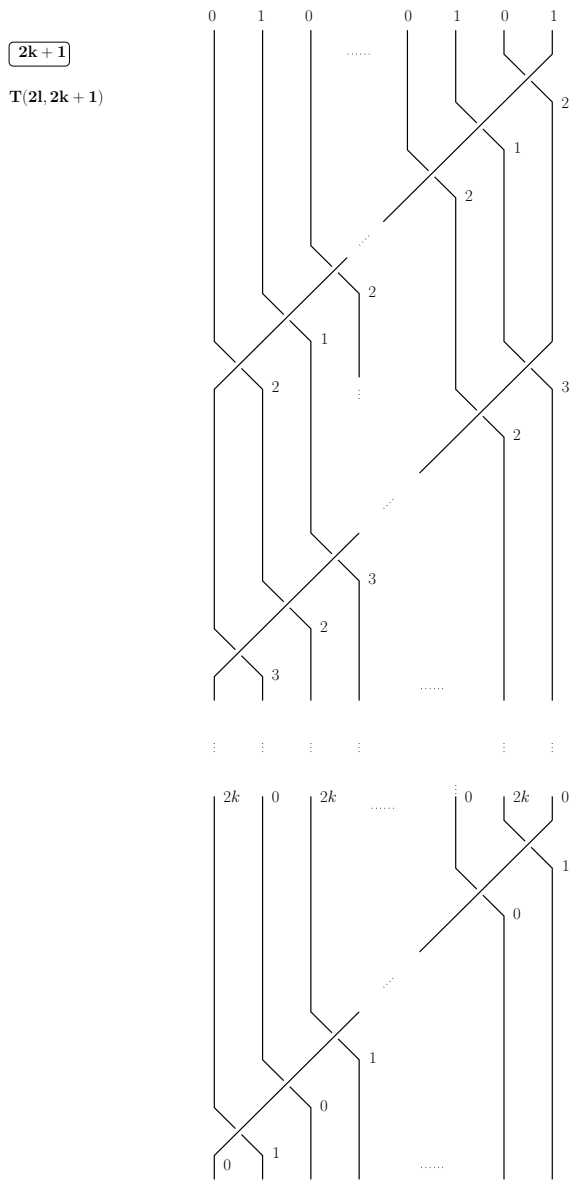}}}
	\caption{A non-trivial $(2k+1)$-coloring of a $T(2l, 2k+1)$ using all colors available: for any positive integer $l$ this Figure stands for a non-trivial $(2k+1)$-coloring.}\label{fig:t2l2k+1}
\end{figure}
\item View Figure \ref{fig:t2k+12l}.  Again, we believe this Figure clarifies the existence of a $(2k+1)$-coloring using all colors available, as stated. Now for the issue of color uniformity. Consider the sectors between two consecutive bridges. In each of these sectors we find the sequence of colors (from left to right) $1, 2, 3, 4, \dots , 2k, 0$ or $0, 2k, 2k-1, \dots , 3, 2, 1$. Now $0$'s and $1$'s occur on arcs belonging to two such consecutive sectors; the remaining colors occur on arcs which belong to only one of these sectors. Thus colors $2, 3, 4, \dots , 2k-1, 2k$ occur each $2l$ times. Colors $0$ and $1$ occur each $l$ times. In the case $\gcd (2k+1, 2l) = 1$ (the knot case) note that with this diagram it is not possible to have a $(2k+1)$-coloring using all colors available and with a uniform distribution of the colors. In fact, this diagram has $2k\times 2l$ crossings (thus arcs) and this number is not divisible by $2k+1$ if $\gcd (2k+1, 2l) = 1$.
\begin{figure}[!ht]
	\psfrag{T2k+12l}{\Huge$\mathbf{T(2k+1, 2l)}$}
	\psfrag{2k+1}{\Huge$\mathbf{2k+1}$}
	\psfrag{2k}{\Huge$2k$}
	\psfrag{2k-1}{\Huge$2k-1$}
	\psfrag{2k-2}{\Huge$2k-2$}
	\psfrag{0}{\Huge$0$}
	\psfrag{1}{\Huge$1$}
	\psfrag{2}{\Huge$2$}
	\psfrag{3}{\Huge$3$}
	\psfrag{4}{\Huge$4$}
	\psfrag{5}{\Huge$\mathbf{5}$}
	\psfrag{6}{\Huge$6$}
	\psfrag{7}{\Huge$\mathbf{7}$}
	\centerline{\scalebox{.25}{\includegraphics{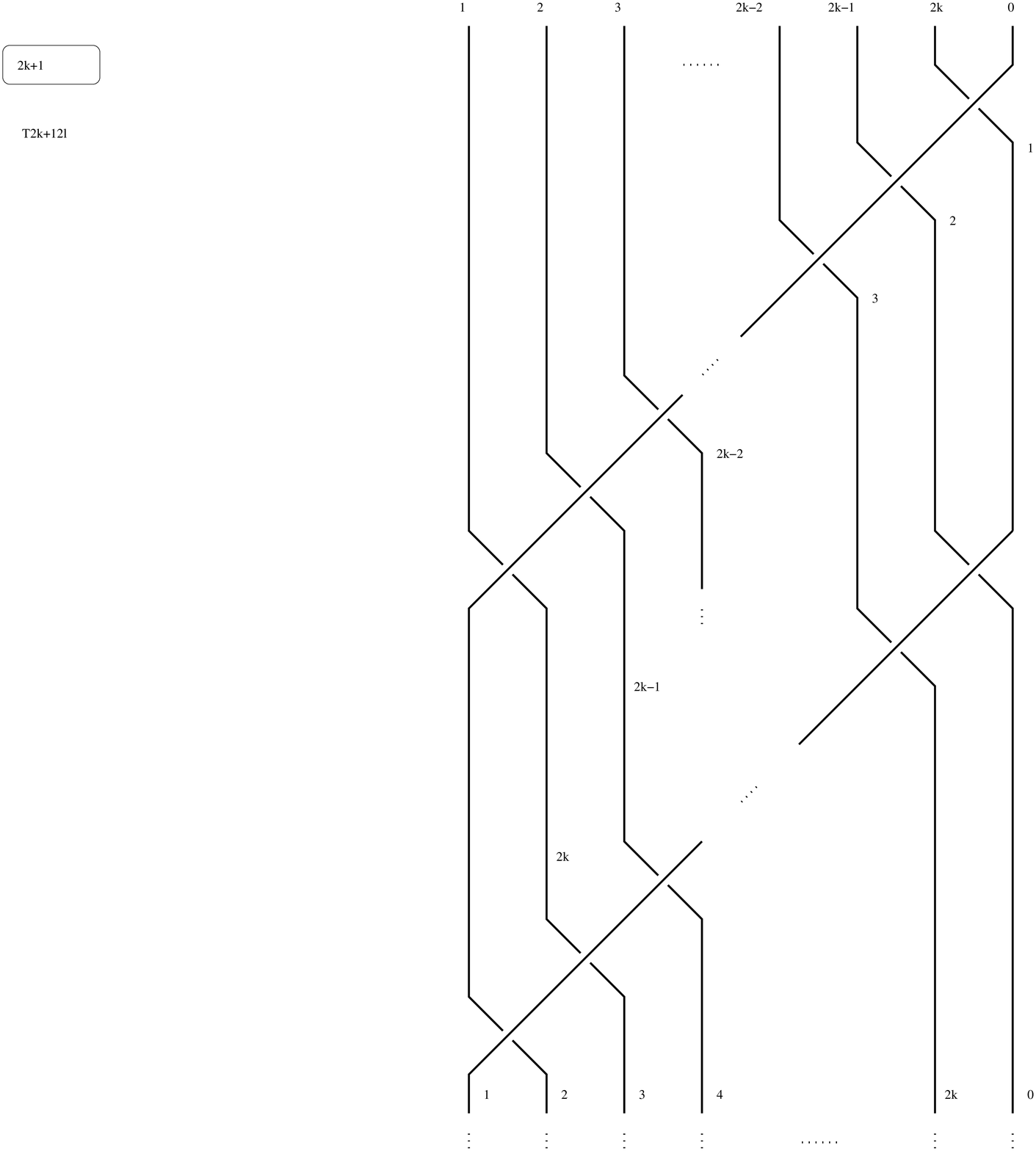}}}
	\caption{A non-trivial $(2k+1)$-coloring of a $T(2k+1, 2l)$ using all colors available: for any positive integer $l$ this Figure stands for a non-trivial $(2k+1)$-coloring.}\label{fig:t2k+12l}
\end{figure}
\item Suppose $k>1$. There are $(2k-1)\times 2l$ crossings in the diagram. Since $2kl$ does not divide $(2k-1)\times 2l$, non-trivial $2kl$-colorings using all colors available cannot exhibit a uniform distribution.

    If $k=1$, then we are dealing with torus knots of type $(2, 2l)$ for which a color input $(0, 1)$ on the top of the diagram will propagate downwards producing a $2kl$-coloring of $T(2k, 2l)$ exhibiting all colors (see \cite{kl}, for instance).
\end{enumerate}
$\hfill \blacksquare$

\bigbreak

We leave the following question for future consideration.
If $k>1$ is there a $2kl$-coloring of the closure of the $(2k)$-braid $(\sigma_{2k-1}\sigma_{2k-2} \cdots \sigma_1)^{2l}$ using all colors available?

\bigbreak

\section{Proof of Theorem  \ref{thm:spec} and Corollary \ref{cor:specbis}.} \label{sect:specpL}

\noindent

Proof (of Theorem \ref{thm:spec}).

\bigbreak
We keep the notation and terminology of the statement of the Theorem.

Consider a diagram of $L$ along with a $p$-coloring of it which uses $\mathrm{mincol}_p L$ colors. Without loss of generality this colored diagram includes a crossing as depicted in Figure \ref{fig:x01}. We call this colored diagram $D_0$.

Using the technique introduced in the proof of Theorem \ref{thm:max} (see also Figures \ref{fig:x01}, \ref{fig:x01II}, and \ref{fig:x01IIn}), we can obtain any other color we need even if along the way we have to obtain repeats first.

For instance, if $3$ is not a color in $D_0$, we can obtain a new diagram from it with color $3$ using the transformation that takes from Figure \ref{fig:x01} to Figure \ref{fig:x01II}. Having performed this transformation even in the case where $3$ is already a color in $D_0$, the stage is set for the introduction of color $4$ via the indicated technique. This procedure allows us to obtain, from $D_0$, a diagram with any number of colors, from $\mathrm{mincol}_p L$ to $\mathrm{maxcol}_p L$. This completes the proof of Theorem \ref{thm:spec}.
$\hfill \blacksquare$
\bigbreak

Proof (of Corollary \ref{cor:specbis}).

This is similar to Corollary \ref{cor:max}. Hypothesis $(1)$ or $(2)$ either ensure there is a non-trivial coloring with a crossing equipped with colors as in Figure \ref{fig:x01}, or provide sufficient conditions for the production of such a coloring. As of this point Theorem \ref{thm:spec} can be applied to prove Corollary \ref{cor:specbis}.$\hfill  \blacksquare$
\bigbreak

\begin{figure}[!ht]
	\psfrag{5}{\Huge$\mathbf{5}$}
	\psfrag{0}{\Huge$0$}
	\psfrag{1}{\Huge$1$}
	\psfrag{2}{\Huge$2$}
	\psfrag{3}{\Huge$3$}
	\psfrag{4}{\Huge$4$}
	\centerline{\scalebox{.35}{\includegraphics{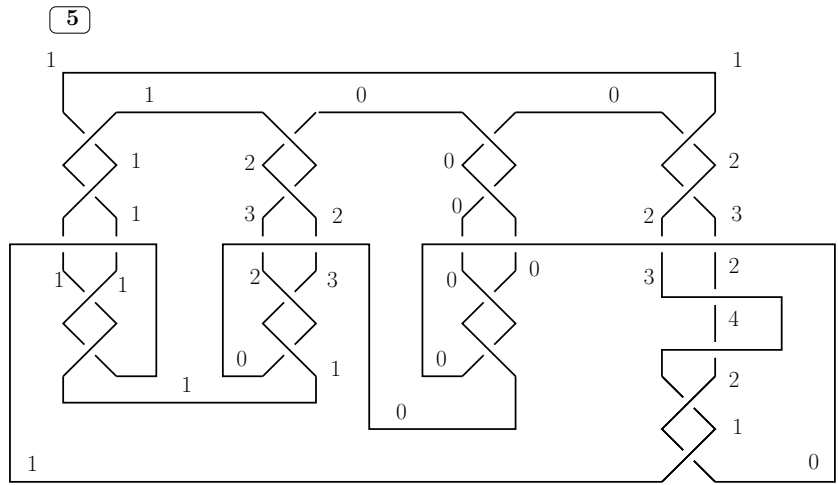}}}
	\caption{Adding color $4$ to a minimal diagram using the technique leaning on a type II Reidemeister move. Total number of colors: $5$.}\label{fig:counter}
\end{figure}

\begin{figure}[!ht]
	\psfrag{5}{\Huge$\mathbf{5}$}
	\psfrag{0}{\Huge$0$}
	\psfrag{1}{\Huge$1$}
	\psfrag{2}{\Huge$2$}
	\psfrag{3}{\Huge$3$}
	\psfrag{4}{\Huge$4$}
	\centerline{\scalebox{.35}{\includegraphics{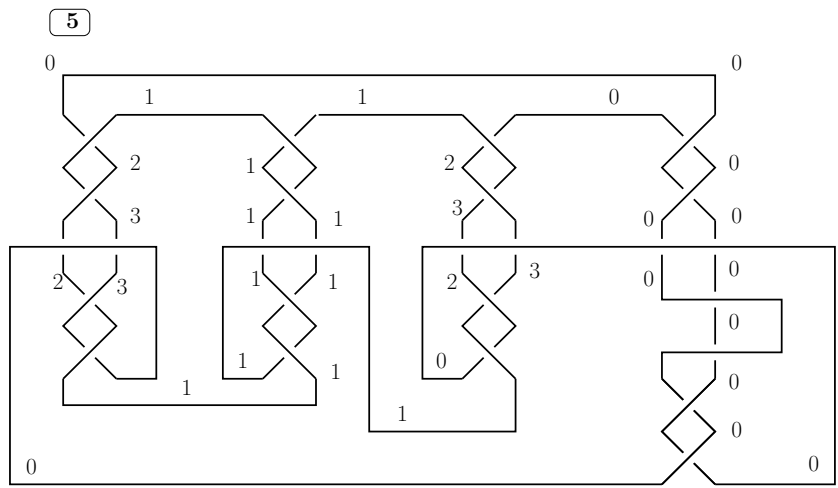}}}
	\caption{The same diagram as in Figure \ref{fig:counter}. Total number of colors: $4$.}\label{fig:counter+}
\end{figure}

We note that the intermediate diagrams in the proof of Theorem \ref{thm:spec} can be equipped with non-trivial colorings with fewer colors than those referred to in the proof, see Figures \ref{fig:counter} and \ref{fig:counter+}.

\bigbreak

\section{Final Remarks.} \label{sect:final}

\noindent

The Kauffman-Harary conjecture on alternating knots of
prime determinant was set forth in \cite{Frank} and proven to be true in \cite{msolis}. Given an odd prime $p$ and an alternating
knot of determinant $p$, any non-trivial $p$-coloring of a reduced
alternating diagram of this knot will assign different colors to different
arcs. Let us call this the KH (Kauffman-Harary) property: a diagram $D$
endowed with a coloring which
assigns different colors to different arcs is said to have the KH property. We will say the diagram satisfies the KH property mod $p$ when it is not clear from context which modulus is meant or to stress the modulus under study.

With this terminology the conjecture above reads: given an odd prime $p$
along with an alternating knot of determinant $p$, any reduced alternating
diagram of this knot possesses the KH property, mod $p$.

From this conjecture we see that from a reduced alternating diagram endowed with a non-trivial coloring modulo an odd prime, we can, in principle, evolve either to maximize or to minimize the number of colors. In this respect we pose the following questions.

\begin{enumerate}
\item Which alternating knots of prime determinant already present the minimum number of colors on reduced alternating diagrams non-trivially colored modulo their determinants?
\item Same as above with ``minimum'' replaced for ``maximum''.  We expect
that the complete list consists of all knots with $det(K)=3$ and torus
knots $T(2,2k+1)$ ($k\ge 2$).
\item Concerning non-alternating link diagrams with prime determinant,
which of them (as diagrams) possess the KH property? Concerning
non-alternating knots with prime determinant, which of their
minimal reduced diagrams posses the KH property? Which non-alternating
knots have at least one minimal reduced diagram which has the KH property?
Consider these same questions when the determinant is not necessarily
prime.
\end{enumerate}

We aim to look into these issues in future work.

\bigbreak

\section{Acknowledgements} \label{sect:ackn}

\noindent

S.J. thanks for support of the Serbian Ministry of Science (Grant No. 174012).

P.L. acknowledges support from FCT (Funda\c c\~ao para a Ci\^encia e a Tecnologia), Portugal, through project number PTDC/MAT/101503/2008, ``New Geometry and Topology''. P.L. also thanks the School of Mathematical Sciences at the University of Nottingham for hospitality during his stay there where the present article was written.


\begin{thebibliography}{99}

\bibitem{breil-oesp-taal} A.-L. Breiland et al, \emph{p-coloring classes of torus knots}, pre-print found at \\
\texttt{http://educ.jmu.edu/{\raise.17ex\hbox{$\scriptstyle\sim $}}\hskip2pt taalmala/OJUPKT/breil{\_}\hskip2pt oesp{\_}\hskip2pt taal.pdf}


\bibitem{CFox}
    R. Crowell, R. Fox,  \emph{Introduction to knot theory},
    Dover Publications 2008




\bibitem{Frank} F. Harary, L. Kauffman, \emph{Knots and graphs. I. Arc graphs and colorings}, Adv. in Appl. Math. {\bf 22}
(1999), no. 3, 312-337







\bibitem{Isidroetal} J. M. Isidro et al, \emph{Polynomials for torus links from Chern-Simons gauge theories}, Nuclear Phys. B {\bf 398}
(1993), no. 1, 187-236





\bibitem{SJablan}
        S. Jablan, R. Sazdanovi\'c , \emph{LinKnot-- Knot Theory by Computer},
        Series on Knots and Everything 21,
        World Scientific Publishing Co., River Edge, NJ 2007





\bibitem{lhKauffman}
        L. H. Kauffman, \emph{Knots and physics}, third ed.,
        Series on Knots and Everything 1,
        World Scientific Publishing Co., River Edge, NJ 2001




\bibitem{kl} L. Kauffman, P. Lopes, \emph{On the minimum number of colors for knots}, Adv. in Appl. Math.
{\bf 40} (2008), no. 1, 36-53



\bibitem{kldmtcs} L. Kauffman, P. Lopes, \emph{Determinants of rational  knots}, Discrete Math. Theor. Comput. Sci.
{\bf 11} (2009), no. 2, 111-122


\bibitem{klgame} L. Kauffman, P. Lopes, \emph{The Teneva game}, J. Knot Theory Ramifications, to appear, DOI: 10.1142/S0218216512501258

\bibitem{LinKnot}
        S. Jablan, R. Sazdanovi\'c , \emph{LinKnot}, \texttt{http://www.mi.sanu.ac.rs/vismath/linknot/}


\bibitem{lm}
        P. Lopes, J. Matias, \emph{Minimum number of Fox colors for small primes}, J. Knot Theory Ramifications, {\bf 21} (2012), no. 3, 1250025 (12 pages)


\bibitem{lm2}
        P. Lopes, J. Matias, \emph{Minimum number of colors: the Turk's head knots case study}, arXiv:1002.4722, submitted

\bibitem{msolis}
        T. Mattman, P. Solis, \emph{A proof of the Kauffman-Harary conjecture},  Algebr. Geom. Topol. {\bf 9} (2009), 2027--2039






\bibitem{Oshiro} K. Oshiro, \emph{Any 7-colorable knot can be colored by four colors}, J.  Math. Soc. Japan, {\bf62}, no. 3 (2010), 963--973









\bibitem{Saito} M. Saito, \emph{The minimum number of Fox colors and quandle cocycle invariants},
 J. Knot Theory Ramifications, {\bf 19}, no. 11 (2010),  1449--1456


\bibitem{satoh} S. Satoh, \emph{5-colored knot diagram with four colors}, Osaka J. Math., {\bf46}, no. 4 (2009), 939--948


\bibitem{Schubert} H. Schubert, \emph{Knoten mit zwei Br\"{u}cken}, Math. Zeit., {\bf65},  (1956), 133--170






\end{thebibliography}
\end{document}